\documentclass[12pt]{article}
\usepackage{amsmath}
\usepackage{amssymb}
\usepackage{amscd}
\usepackage[dvips]{graphics}
\usepackage{amsfonts}

\newcommand{\refm [1]} {(\ref{#1})}

\newcounter{theorem}

\newtheorem{proposition}{Proposition}
\newtheorem{theorem}{Theorem}
\newtheorem{lemma}{Lemma}
\newtheorem{corollary}{Corollary}

\newcommand{\la}{\label}

\newcommand{\PP}{\mathbb{P}}

\newcommand{\nicef}{{\bf {\cal F}}}
\newcommand{\en}{\end{equation}}

\def\qed{\mbox{\rule{0.5em}{0.5em}}}

\def\PP{{\Bbb P}}

\def\BR{{\Bbb R}}

\def\BE{{\Bbb E}}

\def\BC{{\Bbb C}}

\def\qed{\mbox{\rule{0.5em}{0.5em}}}

\def\la{\label}

\def\del1{{\delta_n^{(1)}}}

\newcommand{\ignore}[1]{}

\def\BR{{\Bbb R}}

\def\BE{{\Bbb E}}

\def\BC{{\Bbb C}}

\newcommand{\be}{\begin{equation}}
\newcommand{\eu}{\end{equation}}
\newcommand{\ber}{\begin{eqnarray}}
\newcommand{\ena}{\end{eqnarray}}
\newcommand{\nin}{\noindent}
\newcommand{\non}{\nonumber}

\begin{document}
\title{Asymptotic enumeration by Khintchine- Meinardus probabilistic method: Necessary and sufficient conditions for sub  exponential growth.}
 \author{{\bf Boris L. Granovsky}
\thanks{E-mail: mar18aa@techunix.technion.ac.il} \\
Department of Mathematics, Technion-Israel Institute of Technology,\\
Haifa, 32000,Israel.}\maketitle


\begin{abstract}
In this paper we prove the necessity of the main sufficient condition of Meinardus for sub exponential rate of growth of the number of structures, having multiplicative  generating functions of a general form and establish a new necessary and sufficient condition for normal local limit theorem for aforementioned structures.
The latter result allows to encompass in our study structures with sequences of weights having gaps in their support.
\end{abstract}
\maketitle

{\bf Keywords:} Asymptotic enumeration- Generating function- Partitions-

 Local limit theorem.

{\bf Mathematical subject Classification}

 Primary 05A16-Secondary 05A17,60F05
\vskip.25cm

{\bf I. Introduction and Mathematical setting}

The present  work was motivated by papers \cite{frgr} and \cite{GSE}-\cite{ext} coauthored respectively with Gregory Freiman and Dudley Stark, and by the paper \cite{yif}, by Yifan Yang.
Our objective in this paper is the asymptotic behavior, as $n\to \infty,$ of the quantity $c_n$ depicting the number of combinatorial structures of size $n$.
The mathematical setting below is based on \cite{ext}.
We consider throughout combinatorial objects that decompose into sets of simpler objects called irreducibles, or primes or connected components.
Roughly speaking, the present paper proves the necessity of the main   Meinardus' sufficient condition for the sub exponential rate of growth of $c_n, n\to \infty$ and establishes  quite new necessary and sufficient conditions for the normal local limit theorem. The latter allows to encompass in our study structures with sequences of weights having gaps in their support.

The paper contains four sections. In Section I we provide a background for our study and formulate the mathematical setting, in  Section II we state the  results, which are proven in  Section III.
 The final Section IV contains three examples that hint on  perspectives  for future research.
\vskip.25cm

Let $f$ be a
generating function of a nonnegative sequence $\{c_n,\ n\ge 0,\ c_0=1\}$:
\be
f(z)=\sum_{n\ge 0}c_n z^n,\quad \vert z\vert< 1,\la{sac}
\en
with radius of convergence 1 and such that
\be \lim_{z\to 1^-} f(z)=\infty. \la{triv}\eu
The assumption \refm[triv] implies \be\sum_{n\ge 0} c_n=\infty, \la{sumc}\eu which is a necessary condition for $c_n\to \infty,\ n\to \infty,$ the property that features the models considered in this paper.

Our study is restricted to structures (=models) with generating functions $f$ of the following multiplicative form:
 \be
f=\prod_{k\ge 1}S_k\la{sak1}. \eu
It is appropriate to note that the infinite product  in \refm[sak1] often conforms to $q$- series common in number theory (see \cite{andr},\cite{Kot}). In \cite{Kot} an algorithm was suggested for derivation  the asymptotics of $c_n$ in the case of $q$-series,
 under assumption  that factors of the product have
  asymptotics of exponential type.

We assume that the functions $S_k,\ k\ge 1$ in \refm[sak1] have the following Taylor expansions:  \be
S_k (z)=\sum_{j\ge 0}d_k(j)z^{kj}, \ d_k(j)\ge 0,\ j\ge 0,\  k\ge 1. \la{g3}\eu
By virtue of \refm[sac] and \refm[sak1], the radius of convergence of each one of the series for $S_k$ in \refm[g3] should be  $\ge 1.$
In the literature one can find examples of multiplicative combinatorial structures with radius of convergence of $S_k$
ranging from $1$ to $\infty$ (for references see \cite{ABT}).

 The above setting induces a sequence of multiplicative probability measures (=random structures)  $\mu_n, \ n\ge 1$ on the sequence of sets $\Omega_n, \ n\ge 1$ of integer partitions of $n$, such
 that $c_n$ is a partition function of the measure $\mu_n$ (for more details see \cite{ext}).
The multiplicative measures were introduced in the seminal paper of Vershik \cite{V1} in which he investigated a variety of problems related to limit shapes of the measures $\mu_n$ for classical models of statistical mechanics.
 Our subsequent asymptotic analysis of $c_n,$ as $n\to \infty,$ is based on the probabilistic representation of $c_n,$ derived by Khintchine in \cite{Kh} for classical models of statistical mechanics and then extended in \cite{GSE},\cite{GS},\cite{ext} to  the above defined multiplicative models.
 Khintchine's representation of $c_n,\ n\ge 1,$ which is  identity in the free parameter $\delta>0,$  reads as follows:
\be
c_n=e^{n\delta} f_n(e^{-\delta})\PP\left(Z_n=n\right),\quad n\ge
1,
 \la{rep}\eu
where
\be
f_n=\prod_{k=1}^n S_k
\la{fndef}\eu
is the $n$- truncation of the generic generating function  $f$, and
\be Z_n:=\sum_{k=1}^{n}{Y_k},\quad n\ge 1, \la{znr}\eu where  $Y_k$ are independent integer-valued random variables with distributions derived
from \refm[sak1] and \refm[g3] by setting $z=e^{-\delta},\ \delta>0$ :
\be \PP(Y_k=jk)= \frac{d_k(j)\ e^{-\delta kj}}{S_k(e^{-\delta })},
\quad j\ge 0,\quad k\ge 1.\la{Y}\eu
It is clear from \refm[Y] that the representation \refm[rep] is valid if and only if in \refm[g3] the coefficients $d_k(j)\ge 0,\ k\ge 1, j\ge 0.$

 Khintchine's asymptotic method is based on the    choice of the free parameter $\delta$ in the representation \refm[rep] as a solution, denoted  $\delta_n,$ of the following equation (=Khintchine's equation):
 \be \Big(-\log\nicef(\delta)\Big)^\prime_{\delta}= n, \la{dltn}\eu
where  $\nicef(\delta):=f(e^{-\delta}),\ \delta>0$.
Our subsequent  study is devoted to a special class of multiplicative models defined as follows.

{\bf Definition} A multiplicative model \refm[sak1] is called exponential if under $z=e^{-\tau},\
\Re(\tau)=\delta>0$ its generating function $\nicef(\delta)$ has the expansion:
\be
\nicef(\tau)= \exp\left(\sum_{l=0}^r h_l\tau^{-\rho_l}
-A_0\log\tau +\Delta(\tau) \right), \la{nicef1}
\eu
where
\begin{itemize}
\item $r$ is a given integer;

\item  $\Delta(\tau)<\infty, \tau\in {\cal C}$ is the remainder term that admits  expansion into Taylor series, converging in $\tau\in {\cal C};$

\item
 $h_l> 0,\ l=1,\ldots,r\ \text{and}\ h_0,A_0$ are real constants,
while  $\rho_0=0$ and $ 0<\rho_1<\ldots<\rho_r$ are  positive powers.
\end{itemize}

\begin{proposition}
 The Khintchine's equation \refm[dltn]
has a unique solution $\delta=\delta_n>0$ for all\ \ $n$ sufficiently large, where
 \be \delta_n\to 0, \  n\to \infty.\la{dcuku}\eu Moreover, if the model is exponential, then
\be n\delta_ n\to \infty,\ n\to \infty. \la{limit2}\eu
\end{proposition}
{\bf Proof.} \ We have
\be \Big(-\log\nicef(\delta)\Big)^\prime_{\delta}=\frac{\sum_{k=0}^\infty kc_ke^{-k\delta}}{\sum_{k=0}^\infty c_k e^{-k\delta}}.
\la{bdzr}\eu
The  nominator and the denominator of the fraction in \refm[bdzr] tend to $+\infty$, as $\delta\to 0^+,$ by virtue of the assumption \refm[triv], and it is easy to see that the fraction itself also tends to $+\infty,$ as $\delta\to 0^+.$
Differentiating w.r.t. $\delta$ the RHS of \refm[bdzr] and then applying the Cauchy-Shwartz inequality gives
\be  \Big(-\log\nicef(\delta)\Big)^{\prime \prime}_{\delta}<0,\ \text{for \ all } \delta>0,
\la{seco}\eu
from which  we derive that the fraction in \refm[bdzr] decreases in $\delta>0,$ from $+\infty$ to $0$. This
 proves the existence and the uniqueness of the solution $\delta_n$, as well as \refm[dcuku].
For the proof of  \refm[limit2], we use \refm[nicef1] with $\tau=\delta$, and the fact that in \refm[nicef1], $\rho_r>\rho_{r-1}>\ldots>\rho_1>0,$  to rewrite the equation \refm[dltn] in the case of exponential models as
\be rh_r\delta^{-r-1}_n\sim n, \ n\to \infty. \la{euv}\eu
From the above, \refm[limit2] follows immediately. \qed

 Khintchine (\cite{Kh},p.160)  showed that the solution $\delta=\delta_n$ of \refm[dltn] is the point of minimum of  the entropy of the corresponding model of statistical mechanics (see also \cite{GS} for some more details).

{\bf Historical remark.} In the present paper, as well as in \cite{GSE}-\cite{ext} a combination of Khinchine and Meinardus' asymptotic analysis is employed.  It is interesting to understand  the interplay between the two  methods that originated absolutely independently from each other,\ in the $50$-s of the past century.  Khintchine's objective in \cite{Kh} was calculation of mathematical expectations (rather than $c_n$), with respect to the above measure $\mu_n,$  of such quantities common in statistical mechanics, as occupation numbers.  Occupation numbers  depict  numbers of particles that are at a certain energy level $l, \ l\ge 1.$
The expectations of occupation numbers can be expressed as functions of the ratios $\frac{c_{n-l}}{c_n},$ which by virtue of   the representation \refm[rep] do not depend, as $\delta\to 0,\ n\to \infty,$ on the second factor  in \refm[rep]. Because of it, Khinchine did not need the asymptotic analysis of the second factor in \refm[rep]. The latter asymptotic analysis (for the case of weighted partitions)  was  developed by Meinardus (see \cite{andr}) who proposed to use the
 Mellin transform.
From the other hand, Meinardus used complicated  technique of the saddle point method  that was replaced by Khintchine with his elegant
local limit theorem approach.  \qed

As in \cite{ext}, we restrict the study to functions $S_k,\ k\ge 1$ of the specific form
\be\label{frame}
S_k(z)= \big(S(a_k z^k)\big)^{b_k},
\en
where  the series \be S(z)=\sum_{j=0}^{\infty} d_jz^j, \ d_j\ge 0,\ j\ge 0 \la{szs}\eu has a radius of convergence $\ge 1$ and where $0< a_k\le 1$, $ b_k\ge 0,\ k\ge 1$ are given sequences  of the two parameters of the model. Note that  in \refm[szs]  $d_0=1$, by virtue of \refm[sac], \refm[sak1] and the fact that $c_0=1$.

In \cite{ext} it was described the combinatorial meaning of the parameters $a_k,b_k$. In the sequel of this section we mention another interpretation of these parameters.

For multiplicative models, the aforementioned setting is also used  for the study of another asymptotic problem,
which is a limit shape, the topic which has a rich history.
We
 mention below two recent papers on limit shapes.
Developing the work \cite{V1}, Yakubovich (\cite{YA}) derived  limit shapes for models \refm[frame] in the case $a_k=1,\ k\ge 1$, under some analytic conditions on the function $S$ and on the parameters $b_k,\ k\ge 1$, while Bogachev in \cite{bog}
developed a unified approach to derivation of limit shapes in the case of equiweighted parts: $b_k=b> 0,\ k\ge 1.$
In the sequel  of the present paper, some  other links to research on limit shapes will be indicated.

To simplify the exposition, we make the additional assumption on $a_k:$
\be a_k>0,\ k\ge 1. \la{aks}\eu
We note that the assumption  $d_j\ge 0,\ j\ge 0$ in \refm[szs]
is not sufficient for $c_n\ge 0, \ n\ge 1$. In fact, in the case of weighted partitions  with distinct parts, we have $S(z)=1+z,$ so that  not all   coefficients of the binomial series $\big(S(a_kz^k)\big)^{b_k}, \ k\ge 1$ are nonnegative, unless $b_k,\ k\ge 1$ are integers. In particular, in the case $a_k=1,\ k\ge 1$, it is not difficult to check that $c_2=\frac{b_1(b_1-1)}{2}+ b_2<0,$ for some values of $0<b_1<1,\ b_2>0.$
In this connection recall that for unrestricted weighted partitions, the property $c_n\ge 0,\ n\ge 1$ holds for all $b_k\ge 0,\ k\ge 1.$

  By  the  assumptions made,
 $\log S(z)$ can be expanded as
\be\label{logSexp}
\log S(z)=\sum_{j=1}^\infty \xi_j z^j,
\en
with the  radius of convergence $\ge 1.$ In view of \refm[sac], the function $S(z),$ in the case of  an exponential model, may have zeros and singularities on the unit circle $\vert z\vert=1$ only. We also point  that $ \xi_j\ge 0,\ j\ge 1$ is a necessary and sufficient condition for $c_n\ge 0,\ n\ge 0$ to hold under all $b_k\ge 0,\ k\ge 1,$ which explains the aforementioned dichotomy between unrestricted partitions and partitions into distinct parts.

We assume further on that all  singular points of $S(z)$, if they exist, are poles $z_0:\vert z_0\vert =1,$
which means that
\be
 S(z)\sim \frac{L(z)}{( z-z_0)^l},\ z\to z_0: \vert z\vert<1,
\la{rgul}\eu
with a given integer $l\ge 1$ and with a function $L$  analytic  in the unit disk and such that $0< \vert L(z)\vert<\infty,\ \vert z\vert\le 1 .$
The assumption \refm[rgul] conforms to the one by Yakubovich in \cite{YA}, in the particular case $z_0=1$
and $L$ is a slowly varying function.  Assumption \refm[rgul] extends our study to models with more general $S(z)$, e.g.  $S(z)=\frac{1+z}{1-z^p}$ with an integer $p\ge 1$ and $\vert z\vert\le 1$.

In connection with \refm[rgul], it is appropriate to recall Vivanti-Pringsheim theorem (see e.g. \cite{hil}) which says that
if the radius of convergence of the series \refm[szs] for $S(z)$ equals to $1,$ then  the assumption $d_j\ge 0,\ j\ge 0$ implies  that
$z=1$ is a singular point of the series \refm[szs].

{\bf Remark}  The  example $S(z)=\exp(\frac{1}{1-z}),\ \vert z\vert<1$ demonstrates that a non regular growth of $S(z),\  z\to 1^-$ may lead to a non exponential growth of $c_n,\ n\to \infty.$
\qed

Setting $\tilde{b}_k=lb_k$ and $\tilde{L}(z)= (L(z))^{1/l}$ in  \refm[frame] we
may assume without loss of generality that  $l=1$ in \refm[rgul].

By \refm[frame] and \refm[logSexp], the following expansion of $\log f(z)$ is valid:
\be
\log f(z) = \sum_{k\ge 1}b_k \log S(a_kz^k):=\sum_{k\ge 1} \Lambda_k z^k,\quad \vert z\vert<1,\la{lofl}\eu
with
\be \Lambda_k=\sum_{j\mid k}b_j a_j^{k/j} \xi_{k/j}.\label{logflam}
\en
{\bf Remark}\ \  In view of \refm[lofl],   the function $\log f(z)$ has a pattern of a harmonic sum with  base functions
$\log S(a_kz^k),\ k\ge 1,$ amplitudes $b_k, \ k\ge 1$ and frequencies $a_k,\ k\ge 1.$
Harmonic sums are widely applied in computer science. For more details see \cite{flaj} which studies the asymptotics of harmonic sums with the help of Mellin transform.

By virtue of  \refm[logflam], the Dirichlet generating function $D$ for the sequence $\Lambda_k,\ k\ge 1$
\be\label{Dirdef}
D(s)=\sum_{k=1}^\infty \Lambda_k k^{-s},
\en
  can be written as the double Dirichlet series:
\be \label{dpres}
D(s)=\sum_{k=1}^\infty\sum_{j=1}^\infty b_k \xi_j a_k^j
(jk)^{-s}.
\en
 It is known that if  $D(s)$ converges  in the half-plane  $\Re( s)>\rho_r,$ for some  $\rho_r>0,$  then $D(s)$ converges absolutely
for $\Re( s)>\rho^*, \rho_r<\rho^*< \rho_r+1.$ This amounts to say that $\rho_r$ is the rightmost pole of the Dirichlet series
$D(s)$ and that $\Lambda_k=o(k^{\rho_r}),\ k\to \infty.$

In the particular case,
 $a_k\equiv a,\ 0<a\le 1 $  the function $D(s)$ can be factored as
\be\label{factor}
D(s)=D_b(s)D_{(\xi,a)}(s),
\en
where
$$
D_{(\xi,a)}(s)=\sum_{j=1}^\infty a^j\xi_j j^{-s}
$$
and
\be\label{Dbdef}
D_b(s)=\sum_{k=1}^\infty b_k k^{-s}.
\en

In the general case,  consider the Dirichlet generating function for the sequence $\{\xi_ja_k^j,\ j\ge 1, \ k\ \text{is fixed}\}:$
\be D_{(\xi,a_k)}(s):=\sum_{j\ge 1}\frac{\xi_ja_k^j}{j^s},\quad 0< a_k\le 1,\quad k\ge 1, \la{dirich}\eu
which allows to rewrite \refm[dpres] as \be D(s)=\sum_{k\ge 1} \frac{b_k}{k^s}D_{(\xi,a_k)}(s). \la{Dsd}\eu

We will assume throughout
 the rest of the paper that the function $D(s),\ s=\sigma+it$ is
 of finite order in the whole domain of its definition, which means (see e.g. \cite{titch}) that
 \be  D(s)= O(\vert t\vert^C),\ t\to \infty,\la{xsr}\eu for some constant $C>0$,  uniformly for $\sigma$ from  the domain of definition of $D(s)$. (Here and throughout the paper $$f(x)= O(g(x)),\ x\to a\in \BR,\ g(x)>0 $$ means that $\vert f(x)\vert \le C_1g(x),$ with some constant $C_1>0$ and for all $x$ sufficiently close to $a$).

Recalling that a Dirichlet series is of finite order in any half-plane from the half-plane of its convergence,  the assumption \refm[xsr] requires that if the Dirichlet series \refm[Dsd] admits analytic continuation, then \refm[xsr] holds also in the extended domain.
Finally, recall that the  assumption \refm[xsr] appears in Meinardus' theorem (see \cite{andr}) as one of the sufficient conditions for  sub exponential growth of $c_n$.
\vspace{0.5cm}

{\bf II. Two main theorems}
\vskip .5cm
\begin{theorem}
A multiplicative model with functions $S_k,\ k\ge 1$ of the form \refm[frame] is exponential if and only if the following two conditions hold:

\begin{itemize}

\item {\bf Condition $I$}
The Dirichlet series $D(s)=\sum_{k=1}^\infty \Lambda_k k^{-s}, s>\rho_r$
admits  meromorphic  continuation to $\cal C,$ where it is analytic except $r\ge 1$
simple poles $0<\rho_1<\ldots<\rho_r$ with respective residues $ A_1>0,\ldots, A_r>0,$ and may be  a simple pole at $s=0$ with residue $A_0.$
  The Taylor expansion of the remainder term $\Delta(\tau)$ in \refm[nicef1] is given by
   \be \Delta(\tau)=\sum_{l\ge 1}\frac{(-1)^lD(-l)}{l!}\ \tau^l,\ \tau\in {\cal C}.\la{atrh,}\eu

\item {\bf Condition $II$}


 All\ \ $r\ge 1$  simple positive poles $\rho_1,\ldots\rho_r$ in Condition $I$ belong to the  Dirichlet series $D_b$ in \refm[Dbdef],
while the Dirichlet series $D_{(\xi,1)}$ in \refm[dirich] may have only one simple pole at $0$.
\end{itemize}
 \end{theorem}
 To formulate Theorem 2 below we need the following
\vskip .25cm
\noindent {\bf Definition}
 Let the random variable $Z_n$ be defined as in \refm[znr], \refm[Y], with $\delta=\delta_n$ given by \refm[dltn]. Then we say that for $Z_n$ the normal local limit theorem ($NLLT$) is in force if

 \be
\PP\left(Z_n=n\right)
 \sim\frac{1}{\sqrt{2\pi {\rm Var(Z_n) }}}\ ,\ n\to \infty. \la{locl}\eu

\begin{theorem}
For an exponential model, the NLLT \refm[locl] holds if and only if
the following two conditions on coefficients $d_j,\ j\ge 0$ in \refm[szs] and weights $b_k,\  k\ge 1$ in \refm[frame]
 are satisfied:
 \be gcd\{j\ge 1:d_j>0\}=1\la{,nur}\eu
 and for any integer $q\ge 2$ and $n\to \infty$,\be  \sum_{1\le k\le n,q\!\!\not| k} b_k\ge \left\{
  \begin{array}{ll}
    C\log n,\ C>0, & \hbox{in  Case $(A)$;} \\
   C\log^2 n,\ C>0, & \hbox{in  Case $(B)$,}
  \end{array}
\right. \la{nscon}
\eu
where the cases $(A)$,$(B)$ are as defined below, in the course of proof (see \refm[n,fho] and the Remark after it).

\end{theorem}

\begin{corollary}
Let all conditions of Theorem 2 hold. Then the following asymptotic formula for $c_n$ is valid:
\be c_n\sim \frac{\delta_n^{\frac{\rho_r}{2}+1}}{\sqrt{2\pi}}\exp\left(\sum_{l=0}^r h_l\delta_n^{-\rho_l}
-A_0\log\delta_n + \Delta(\delta_n) + n\delta_n\right), \ n\to\infty,\la{expgr}\eu
where
$\delta_n$ is the unique solution of the Khintchine's equation \refm[dltn], while
  $\Delta(\tau),$ as well as the constants are as defined in Condition $I$.
\end{corollary}

It is important to note that since $\delta_n=O(n^{-\frac{1}{\rho_r+1}}),\ n\to \infty,$ by virtue of \refm[euv], the sub exponential rate of growth of $c_n,$ as determined by  \refm[expgr],  is  \be O(n^{-\frac{\rho_r+2}{2\rho_r+1}})\exp(O(n^{\frac{\rho_r}{\rho_r+1}})),\ n\to \infty,\la{eruc}\eu
whereas the assumption \refm[sac]  requires $c_n=o(e^n),\ n\to \infty.$
Thus, for $\rho_r>0$ sufficiently large, the rate \refm[eruc] of sub exponential growth of $c_n, \ n\to\infty$
approaches  the maximal possible one.

 In this connection
 note that all models of harmonic sums treated in \cite{flaj} exhibit non exponential rate of growth of $c_n$.
\vskip .5cm
{\bf III Proofs}

\begin{itemize}

\item {\bf  Necessity of Condition $I$ of Theorem 1}
\end{itemize}

The key ingredient in our  proof of the necessity of Condition 1 is  the
 forthcoming  Lemma 1. The lemma is  an obvious extension to our setting of Yifan Yang's Lemma 2 in \cite{yif}, where it is formulated for a special case of  partitions into powers of primes.
\begin{lemma}  \ Let the Dirichlet  series \refm[dpres] converge for $\Re(s)>\rho_r>0$ and
 let $\log \nicef(\delta)$ satisfy \refm[lofl],\refm[logflam]. Then
\be \int_0^{1} \delta^{s-1}\log \nicef(\delta)d\delta= \Gamma(s)D(s)- W(s), \quad \Re(s)>\rho_r, \la{Yan}\eu
where $W(s),\ s\in {\cal C}$ is an entire function.
\end{lemma}
{\bf Proof} Substituting \refm[logflam], gives
 $$\int_0^{1} \delta^{s-1}\log \nicef(\delta)d\delta= \Big(\int_0^\infty - \int_{1}^\infty\Big)\delta^{s-1}\log \nicef(\delta)d\delta=$$$$ \Big(\int_0^\infty - \int_{1}^\infty\Big)\delta^{s-1}
 \sum_{k=1}^{\infty} b_k\sum_{j=1}^\infty\xi_ja_k^j e^{-\delta kj} d\delta=$$$$
 \sum_{k,j=1}^{\infty}\Gamma(s)\frac{b_k a_k^j\xi_j}{k^sj^s} -  \sum_{k,j=1}^\infty \Gamma(s,jk)\frac{b_ka_k^j\xi_j}{k^sj^s},$$ \be\Gamma(s,u):=\int_u^\infty x^{s-1}e^{-x}dx,
 \la{relat12} \eu
where the first double series in \refm[relat12] converges to $\Gamma(s)D(s),$ for $\Re(s)>\rho_r$ and converges absolutely for
$\Re(s)>\rho^*,$ where $\rho^*$ is the abscissa of absolute convergence of the  Dirichlet series \refm[dpres].
 We will show that the second double series in \refm[relat12]
defines an entire function in $s\in {\cal C}.$ For this purpose we use the following bound on the incomplete Gamma function $\Gamma(s,u)$ which itself is entire in $s\in {\cal C}$ for all $u>0.$ Letting  $\sigma=\Re(s),$ we have for any reals $\sigma_1<\sigma_2,$
$$\vert \Gamma(s,u)\vert \le \int_u^\infty x^{\sigma-1}e^{-x}dx\le \Big(\max_{x\ge 1}\max_{\sigma\in [\sigma_1,\sigma_2]}x^{\sigma-1}e^{-\frac{x}{2}}\Big)\int_u^{\infty}e^{-\frac{x}{2}} dx  =$$\be 2C(\sigma_1,\sigma_2)e^{-\frac{u}{2}},
\quad u\ge 1,\la{gbound}  \eu uniformly for
 $\sigma\in [\sigma_1,\sigma_2],$
where $0<C(\sigma_1,\sigma_2)<\infty$
denotes the maximum in \refm[gbound].\\
Next, the   absolute convergence in the half-plane $\Re(s)>\rho^*>0,$ of the double series  $\sum_{k\ge 1,j\ge 1}\frac{b_k}{k^s}\frac{\xi_j a_k^j}{j^s} $
that represents the function $D(s)$ in the above  half-plane, implies  $$\frac{b_k}{k^{\rho^*+\epsilon}}\frac{\xi_j a_k^j}{j^{\rho^*+\epsilon}}\to 0,\ k,j\to \infty, \ \epsilon>0.$$ As a result, applying \refm[gbound] we get
$$\vert W(s)\vert=\vert\sum_{k,j=1}^\infty \Gamma(s,jk)\frac{b_ka_k^j\xi_j}{k^sj^s}\vert\le 2C(\sigma_1,\sigma_2)\sum_{k,j=1}^\infty \frac{b_ka_k^j\vert\xi_j\vert}{(kj)^{\sigma}} e^{-\frac{kj}{2}}=$$$$ 2C(\sigma_1,\sigma_2)\sum_{k,j=1}^\infty \frac{o((kj)^{\rho^*+\epsilon})}{(jk)^{\sigma}}\ e^{-\frac{kj}{2}}<\infty,$$
uniformly for   $\Re(s)=\sigma\in [\sigma_1,\sigma_2],\ \text{with\ any \ reals}\ \sigma_1,\sigma_2$. This proves that the function $W(s)$ is entire.\ \ \ \qed

Assuming that the structure is exponential, it follows from the definition
\refm[nicef1]  with $\tau=\delta_n>0,$  the necessity of the asymptotic formula
$$ \nicef(\delta_n)= \exp\left(\sum_{l=0}^r h_l\delta_n^{-\rho_l}
-A_0\log(\delta_n)+\Delta(\delta_n) \right),$$ \be  \Delta(\delta_n)\to 0,\ n\to\infty. \la{niceff}\eu
By Lemma 1, \be D(s)= \frac{1}{\Gamma(s)}\Big(\int_0^{1} \delta^{s-1}\log \nicef(\delta)d\delta + W(s)\Big),\quad \Re(s)> \rho_r. \la{intgr}\eu Next, substituting  \refm[niceff]  and the Taylor expansion of $\Delta(\delta_n),$  into the integral in \refm[intgr],
 we obtain
\be D(s)= \frac{1}{\Gamma(s)}\Big(\frac{h_0}{s}+ \frac{A_0}{s^2} + \sum_{l=1}^r \frac{h_l}{s-\rho_l} + W(s)+ \sum_{k\ge 1}\frac{\Delta^{(k)}(0)}{k!(s+k)}\Big),\ s>\rho_r. \la{DSA}\eu
Since $\frac{1}{\Gamma(s)}$ is an entire function, with zeros at $s=-n, n\ge 1$, \refm[DSA] says that  the function $D$ is analytic in $\BC$, except the simple positive poles $\rho_1,\ldots, \rho_r,$ with the respective residues  $ A_l=\frac{h_l}{\Gamma(\rho_l)}> 0, \ l=1,\ldots,r$ and a simple pole at $s=0$ if  $A_0\neq 0$ in \refm[nicef1].  The latter, together  with the fact that
$s=-n,\ n\ge 1 $ are simple poles of $\Gamma(s)$ with residues $\frac{(-1)^n}{n!},$ respectively, gives $$\Delta^{(k)}(0)=(-1)^kD(-k),\ k\ge 1.$$
As a result, the proof of the necessity of Condition $I$ is completed.

\begin{itemize}
\item {\bf  Necessity of Condition $II$ of Theorem 1}
\end{itemize}

Our  first objective  is to prove the remarkable fact that $\rho_1,\ldots,\rho_r$ are poles of the Dirichlet generating series
$D_b(s)=\sum_{k\ge 1}\frac{b_k}{k^s}$ for the weights $b_k,\ k\ge 1.$

For this purpose we firstly prove that for a given $k\ge 1$ and  $0<a_k\le 1,$ the Dirichlet series $D_{(\xi,a_k)}$ given by \refm[dirich]  has no positive poles.  The two cases $0<a_k<1$ and $a_k=
 1$ should be distinguished.

In the first case, $D_{(\xi,a_k)}(0)=\log S(a_k)<\infty,\ 0<a_k<1,$ since  the radius of convergence of the series \refm[logSexp] is $\ge 1$. Consequently, the Dirichlet series $D_{(\xi,a_k)}(s),\ 0<a_k<1$ converges in the half- plane $\Re(s)\ge 0,$ and therefore it is analytic in this domain.

If  $a_k=1$ for a given $k,$ there are the  following two possibilities:

$(i)$ $D_{(\xi,1)}(0)= \sum_{j\ge 1}\xi_j<\infty.$ This says that in the case considered the Dirichlet series $D_{(\xi,1)}(s)$ is analytic in the half-plane $\Re(s)\ge 0.$

{\bf Example} Partitions into distinct parts: $$S(z)=1+z,\ D_{(\xi,1)}(0)=\sum_{j\ge 1}\xi_j= \sum_{j\ge 1} \frac{(-1)^{j-1}}{j}<\infty.$$

$(ii)$ $ \sum_{j\ge 1}\xi_j= \infty,$
 which amounts to saying that  $z=1$ is the singular point of the series \refm[logSexp], with radius of convergence $1$ and with $\log S(z)\vert_{z\to 1^-}\to \infty$.
Consequently, $S(z)\to \infty,\quad z\to 1^-,$  so that the assumption \refm[rgul] is in force. Recalling that in \refm[rgul] it   can be taken $l=1$,  we have
 \be \log S(z)=\sum_{j\ge 1} \xi_j z^j \sim \log L(z)+ \log(\frac{1}{1-z})\sim \log(\frac{1}{1-z}),\  z\to 1^-,  \la{regl3}\eu
where the second $"\sim" $ is by   the  properties of the  function $L$ as stated in \refm[rgul].
Thus, we obtain from \refm[regl3], \be \big(\log S(z)\big)_z^\prime=\sum_{j\ge 1} j\xi_j z^{j-1}\sim \frac{1}{1-z},\quad z\to 1^-. \la{logs}\eu
We apply now Karamata's tauberian theorem (see e.g. \cite{fel}) to the  asymptotic relation  in \refm[logs] to derive
that $j\xi_j\sim 1, \ j\to \infty.$ Thus, $\xi_j\sim \frac{1}{j}, \ j\to \infty,$
which implies that in the case considered:
  \be D_{(\xi,1)}(\rho)=\sum_{j\ge 1}\xi_jj^{-\rho}<\infty, \quad \text{for all}\ \rho>0,
 \la{slow}\eu while $D_{(\xi,1)}(0)=\infty.$
This says that in case $(ii)$,  under the assumption \refm[rgul] (with $l=1$), the Dirichlet series $D_{(\xi,1)}(s)$ has in the half-plane $\Re(s)\ge 0$ only one simple   pole at $s=0.$

{\bf Example} Unrestricted partitions: $$S(z)=\frac{1}{1-z},\ \ D_{(\xi,1)}(0)= \sum_{j\ge 1} \frac{1}{j}=\infty.$$
The corresponding
Dirichlet series $D_{(\xi,1)}(s)=\zeta(1+s)$ has a unique simple pole at $s=0.$

 From the above proven fact that for a given $k\ge 1$ the function $D_{(\xi,a_k)}$ has no positive poles we will  derive now that all positive poles of $D(s)$ belong to the Dirichlet series $D_b(s)$.
 If $a_k\equiv 1$, then the claim follows immediately from \refm[factor]. In the general case, $D_{(\xi,a_k)}(\rho)<\infty,\ k\ge 1, \rho>0$ implies $\sup_{k\ge 1}D_{(\xi,a_k)}(\rho):=u(\rho)<\infty,\ \rho>0,$ because $0<
  a_k\le 1$. In view of this, taking $\rho=\rho_l>0,\ 1\le l\le r$ we have
 $$\infty=D(\rho_l)\le u(\rho_l)\sum_{k\ge 1}\frac{b_k}{k^{\rho_l}},$$
which says that $\rho_l,\ l=1,\ldots, r$ are indeed the poles of $D_b(s)$.\qed

{\bf Remarks}
$(i)$
 The property of $D_{(\xi,a_k)}$ stated in Condition $II$ is shared
by the three classic combinatorial structures, which are multisets, selections and assemblies (see \cite{ABT},\cite{GSE},\cite{GS}).

$(ii)$ A Dirichlet series with real coefficients may have complex poles. For example,
$$\frac{\zeta(s+1)}{\zeta(s)}=\zeta(s+1)\frac{1}{\zeta(s)}=\zeta(s+1)\sum_{k=1}^\infty \frac{\mu(k)}{k^s},\quad \Re(s)>1, $$
where $\mu(k)$ is the Moebius function. The LHS of the above relation is  a product of two Dirichlet series $\zeta(s+1)$ and $1/\zeta(s)$ with real coefficients. However, the product   has complex poles which are complex zeros of $\zeta(s)$ on the critical line $s=1/2+it.$

\begin{itemize}
\item {\bf   Sufficiency of Conditions $I$ and $II$ of Theorem 1}
\end{itemize}

 Our proof of sufficiency of Condition $I$ for \refm[nicef1] follows  Meinardus' scheme which is based on application of Mellin transform. Here we sketch the scheme, assuming that the  details can be found in  \cite{ext}.
We use the fact that
$e^{-u}$, $\Re(u)>0$, is the Mellin transform of the Gamma function:
\be e^{-u}=\frac{1}{2\pi i}\int_{v-i\infty}^{v+i\infty}
u^{-s}\Gamma(s)\,ds,\quad  \Re(u)>0, v>0. \la{Mellin} \en
Applying \refm[Mellin] with $u=\tau: \Re(\tau)=\delta>0$ and $v=\rho_r+\epsilon,\ \epsilon>0$  we have\begin{eqnarray}
\log
~\nicef(\tau)
&=&
\sum_{k=1}^\infty b_k\log S\left(a_k e^{-\tau k}\right)\nonumber\\
&=&
\frac{1}{2\pi i}\int_{\epsilon+\rho_r-i\infty}^{\epsilon+\rho_r+i\infty} \tau^{-s}\Gamma(s)D(s)ds, \la{intrep3}
\end{eqnarray}
where $D(s)$ is a meromorphic  continuation to ${\cal C}$ of the Dirichlet series \refm[dpres].

Next, assuming that the Condition $I$ holds and recalling \refm[xsr], we apply the  residue theorem  for the integral \refm[intrep3] in the complex domain $\Re( s)\le \rho_r+\epsilon,$ to get the formula:
 \be \log~\nicef(\tau)= \sum_{l=0}^r h_l\tau^{-\rho_l}
-A_0\log\tau +\Delta(\tau), \la{tsdr}\eu
where $h_l=A_l\Gamma(\rho_l), \ l=1,\ldots,r$ and where the expansion \refm[atrh,] of the remainder term $\Delta(\tau)$ follows from the fact  that in the domain $\Re( s) <0$ the integrand  $\delta^{-s}D(s)\Gamma(s)$ has simple poles at $s=-k,\ k=1,2,\ldots,$ only.
Exponentiating \refm[tsdr] gives \refm[nicef1].

{\bf Remark} In the previous research, started from the aforementioned seminal paper by Meinardus it was always assumed that $D$ admits meromorphic continuation to $-C_0<\Re( s)<\rho_1,$ for some $0<C_0<1.$ Accordingly, the error of the asymptotic expansion of
$\log~\nicef(\delta)$ was $O(\delta^{C_0}),$ $\delta\to 0^+.$
\vskip.25cm
The sufficiency of Condition $II$ is obvious.\qed

{\bf Two remarks regarding} \refm[tsdr].\\ \hskip 1cm $(i)$ By  \refm[tsdr] with $\tau=\delta>0,$
\be \frac{\log \nicef(\delta)}{\delta^{-\rho_r}}\to h_r=A_r\Gamma(\rho_r),\ \delta\to 0^+, \eu
since $\rho_r>0$ is the rightmost pole. In the particular case of ordinary partitions the above asymptotics has  deep meanings in statistical physics and combinatorics, being related respectively, to  the shape of a crystal at equilibrium and the limit shape of a random partition of large $n$. In both cases the interpretation is based on treating a limit shape as a solution of a certain variational problem (see \cite{ok} and references therein).

$(ii)$ We consider here the formula \refm[tsdr] in the special case of
   ordinary partitions: $D(s)=\zeta(s)\zeta(s+1),$  $r=1,\ \rho_1=1,\ h_1= \zeta(2) , A_0=-\zeta(0)$.\ \ In the case considered  \refm[tsdr] becomes
   :
\begin{eqnarray}
\log ~\nicef(\delta)
&=&\delta^{-1}\zeta(2)-\zeta(0)\log \delta + \zeta^\prime(0)+\nonumber\\
&& \sum_{l=1}^N \delta^l\frac{(-1)^l}{l!}\zeta(-l+1)\zeta(-l)+ \Delta I(C_N;\delta),
\la{nrs}
\end{eqnarray}
where we denoted
\be \Delta I(C_N;\delta)=\frac{1}{2\pi i}\int_{-C_N - i\infty}^{-C_N + i\infty}\delta^{-s}\Gamma(s)\zeta(s)\zeta(1+s)ds , \la{anj}\en
for  a fixed  integer $N$ and  $C_N=N+ \frac{1}{2}$.
 Next, from the functional equation for zeta - function and from the Gauss multiplication formula for the Gamma function
(see Thm. 15.1 in \cite{AAR}) we have:
$$\Gamma(s)\zeta(s)\zeta(s+1)=(2\pi)^{2s}\zeta(1-s)\zeta(-s)\Gamma(-s),$$
which gives
\be \Delta I(C_N;\delta)=\frac{1}{2\pi i}\int_{-C_N - i\infty}^{-C_N + i\infty}\delta^{-s}(2\pi)^{2s}\zeta(1-s)\zeta(-s)\Gamma(-s)ds . \la{anj}\en
Now we use the asymptotic bound
$$\vert\zeta(u)\vert\le 1+ \sum_{n\ge 2} n^{-v}\le 1 + \int_2^\infty x^{-v}dx=1+ \frac{2^{1-v}}{v-1}=1+ O(2^{-v}v^{-1}),$$ where we denoted $$ 1<v:=\Re(u)\to +\infty.
  $$
This and    the Mellin transform formula of the Gamma function allow to bound $\Delta I(C_N;\delta)$ in \refm[anj]:
\be \vert\Delta I(C_N;\delta)\vert\le e^{-(2\pi)^2\delta^{-1}}(1 + O(2^{-C_N}C_N^{-1}))\to e^{-(2\pi)^2\delta^{-1}}, \ \text{as}\ 0<C_N\to \infty. \non\en
 The latter together with the fact  that $$\zeta(-l+1)\zeta(-l)=0,\ l\ge 1 $$
allows to derive   from \refm[nrs],
\begin{eqnarray} \log ~\nicef(\delta)=
&&\delta^{-1}\zeta(2)-\zeta(0)\log \delta + \zeta^\prime(0)- \zeta(0)\zeta(-1)\delta + \nonumber\\
&& e^{-(2\pi)^2\delta^{-1}},\ \delta\to 0^+. \la{pkd}
\end{eqnarray}
 Exponentiating the last expression  and setting  $x=e^{-\delta}$ recovers the formulae (8.6.1),(8.6.2), p.117 in \cite{Hardy}, for $\nicef(\delta)$.  It is emphasized by Hardy(\cite{Hardy}), that the remainder term $\exp\big(e^{-(2\pi)^2\delta^{-1}}\big)$ in the expansion of  the function $\nicef(\delta)$ goes to $1$ very fast, as $\delta\to 0^+.$ The  formula \refm[pkd], which is the key ingredient of the famous Hardy-Ramanujan expansion for the number of ordinary partitions was derived (see e.g. \cite{Hardy}) in a quite different way, based on the remarkable fact that the generating function for ordinary partitions is an elliptic function obeying  a certain functional equation.  Finally, note that the mysterious exponent $1/24$ in the aforementioned Hardy- Ramanujan  formula is equal to $\zeta(0)\zeta(-1).$

\begin{itemize}
\item {\bf  Theorem 2}
\end{itemize}
  {\bf Three auxiliary facts.}\ \ In $(i)-(iii)$ below we assume that the structure is exponential.  The detailed proofs and the history of $(i)$ and $(ii)$   can be found in \cite{frgr},
  \cite{FrP}. Lemma 2 in $(iii)$  is new.

$(i)$ Asymptotics of $\delta_n,$ as $n\to \infty$.

Firstly we show that in the representation \refm[rep]
with $\delta=\delta_n$,
\be  f_n(e^{-\delta_n})= \nicef(\delta_n)+ \epsilon_n, \ \epsilon_n\to 0, \ n\to \infty. \la{asnicef}\eu
Recalling the expression \refm[logflam] for $\log f(z)$   and that in  \refm[szs], $d_0=1,d_j\ge 0,\ j\ge 1$.
 we denote \be
l_0=\min\{j\ge 1:d_j> 0\}.\la{lopx}
\eu
To avoid the trivial case $S(z)\equiv 1,$ we assume $l_0<\infty.$ Now we have,
 \be S(a_k e^{-\delta_n k})-1= O(a_k^{l_0}e^{-l_0k\delta_n})\to 0,\ n\to \infty,\ \text{for\ all}\ k\ge n,\eu
since $n\delta_n\to \infty$, by Proposition 1.  Consequently,
\be \log S(a_ke^{-k\delta_n})= O(a_k^{l_0}e^{-l_0 k\delta_n })\to 0,\  n\to \infty,\
\text{for\ all}\ k\ge n.\eu
As a result, \be
\sum_{k=n+1}^\infty
 b_k\log  S(a_ke^{-k\delta_n})=
\sum_{k=n+1}^\infty b_k O(a_k^{l_0}e^{-l_0k\delta_n })
 \to 0 ,\quad n\to \infty,\la{aux}\eu
where the last step is because $b_k=o(k^{\rho_r}),\ k\to \infty$, since $\rho_r $ is the rightmost pole of $D_b(s)$  and because of \refm[euv].
\refm[aux] proves \refm[asnicef].
Next, substituting \refm[niceff] into the LHS of the Khitchine's  equation \refm[dltn], produces  the asymptotic expansion of the solution $\delta_n, \ n\to \infty$ of the equation. For the case of multiple poles, i.e. $r>1,$ the
expansion was firstly obtained
in \cite{ext}.
 For our subsequent study  we will need only the main term of the above expansion which is obtained from  \refm[euv]:
  \be \delta_n\sim (\rho_rh_r)^{\frac{1}{\rho_r+1}} n^{-\frac{1}{\rho_r+1}},\ n\to \infty,\la{asdeltan} \eu
where $$ h_r=A_r\Gamma(\rho_r).$$

{\bf Remark} In  connection with \refm[asdeltan] it is in order to note that in the theory of limit shapes the parameter $\delta$ (called there scaling) is taken to be equal $O(n^{-\frac{1}{\rho_r+1}}),$ which is, roughly speaking, \refm[asdeltan] (see e.g \cite{bog},\cite{V1},\cite{YA}).
The aforementioned coincidence is explained by the fact that the derivation of limit shapes  consists of  asymptotic approximation of probabilities with respect to the same multiplicative measure $\mu_n$ as in our setting.

$(ii)$ Representation of $\PP\left(Z_n=n\right).$

  We start from the formula
\be \PP\left(Z_n=n\right)=
 \int_{-1/2}^{1/2}\phi_n(\alpha)
e^{-2\pi in\alpha}d\alpha:=I_1+I_2, \la{nollt}\eu
 where the random variable $Z_n$ is as defined in \refm[znr], \refm[Y] and  $\phi_n(\alpha)$ is the characteristic function of $Z_n,$ while
$$
I_1= I_1(n)=\int_{-\alpha_0}^{\alpha_0}\phi_n(\alpha)e^{-2\pi
in\alpha}d\alpha,\ \text{\ with}\ \alpha_0=\alpha_0(n)=(\delta_n)^{\frac{\rho_r+2}{2}}\log n$$
and
\be
I_2= I_2(n)=\int_{-1/2}^{-\alpha_0}\phi_n(\alpha)e^{-2\pi
in\alpha}d\alpha +\int_{\alpha_0}^{1/2}\phi_n(\alpha)e^{-2\pi
in\alpha}d\alpha. \la{alpha0}
\eu
Our first goal will be to derive  the asymptotics of the integral $I_1=I_1(n),\ \text{as}\ \ n\to\infty.$
Let for a given $n,$  $B_n^2$ and $T_n$ be defined by
\be
Var Z_n:= B_n^2=
\Big(\log f_n\big(e^{-\delta}\big)\Big)^{\prime\prime}_{\delta=\delta_n} \label{BT21}
\end{equation}
and
\be
T_n:=-\Big(\log f_n\big(e^{-\delta}\big)\Big)^{\prime\prime\prime}_{\delta=\delta_n}.\label{BT2}
\end{equation}
Due to the fact that
$$\phi_n(\alpha)e^{-2\pi in\alpha}=\BE \exp{\big(2\pi \alpha i( Z_n-n)\big)},\ \alpha\in \BR,$$
the following   expansion in $\alpha$ is valid, when $n$ is fixed:
 \begin{eqnarray}
 \non\phi_n(\alpha)e^{-2\pi in\alpha}
&=&
\exp{\left(2\pi i\alpha(\BE Z_n-n)-2\pi^2\alpha^2B_n^2+O(\alpha^3)T_n\right)}\nonumber\\
&=&\exp{\left(-2\pi^2\alpha^2B_n^2+O(\alpha^3) T_n\right)},\quad
\alpha\rightarrow0, \non
 \end{eqnarray}
where the second equation is due to \refm[dltn] and the fact that $$\Big(-\log\nicef(\delta)\Big)^\prime_{\delta}=\BE Z_n(\delta),\ \delta>0.$$
 It follows from \refm[niceff]  that  the main terms in the asymptotics for $B_n^2$ and   $T_n$ depend on the rightmost pole $\rho_r>0$ only:
\be
 B_n^2
 \sim
K_2(\delta_n)^{-\rho_r-2},\ n\to \infty, \label{BT21}
\end{equation}
where $K_2= h_r\rho_r(\rho_r+1)$
and \be
T_n\sim K_3(\delta_n)^{-\rho_r-3},\quad n\to \infty, \la{engy}\eu
 where
$K_3=h_r\rho_r(\rho_r+1)(\rho_r+2).$
 Therefore, by the  choice of $\alpha_0$ as in \refm[alpha0],$$B^2_n\alpha_0^2\to \infty, \ T_n\alpha_0^3\to 0, \ n\to \infty.$$
Consequently, by  the same argument as in the proof of the $NLLT$ in  \cite{GSE},
\be\label{I1sim}
I_1\sim \frac{1}{\sqrt{2\pi B_n^2}}\sim (2\pi K_2)^{-1/2}(\delta_n)^{1+\frac{\rho_r}{2}} , \ n\to \infty.
\en
$(iii)$ Bounding  the integral $I_2$, as $n\to \infty$.

 For the $NLLT$ \refm[locl] to hold it is necessary and sufficient that \be I_2=o(I_1), \ n\to \infty.\la{relat}\eu
  \begin{lemma}
For the
 NLLT \refm[locl] to hold it is necessary and sufficient that
\be \vert \phi_n(\alpha)\vert= o(\delta_n^{1+ \frac{\rho_r}{2}}),\ n\to \infty,\ \alpha\in [\delta_n,1/2]. \la{ygb}\eu
\end{lemma}
{\bf Proof}  The following expression is valid for  the  multiplicative models considered:
\be \log\phi_n(\alpha)=\sum_{k=1}^n b_k\Big(\log(S(a_ke^{-\tau_n k})-\log S(a_ke^{-\delta_nk})\Big),\ \ \alpha\in \cal R,\la{ch12}\eu
where $\tau_n=\tau_n(\alpha)=\delta_n-2\pi i \alpha.$
It is easy to see that \refm[asnicef] holds with $\delta_n$ replaced with $\tau_n$.  In view  of \refm[niceff]  we thus have:
$$\log\phi_n(\alpha)=\sum_{k=1}^\infty b_k\Big(\log S(a_ke^{-\tau_n k})-\log S(a_ke^{-\delta_nk})\Big)+\epsilon_n=
\log\nicef(\tau_n)-\log\nicef(\delta_n)+\epsilon_n, $$
\be  \epsilon_n=\epsilon_n(\alpha)\to 0,\  n\to \infty,\ \alpha\in {\cal R},\la{charect}\eu
where for exponential structures
\be \nicef(\tau_n)=\exp\left(\sum_{l=0}^r h_l\tau_n^{-\rho_l}
-A_0\log\tau_n+ \Delta(\tau_n)\right), \la{taudel}\eu
by \refm[nicef1]. Applying  \refm[niceff], \refm[asnicef],\refm[rep] and \refm[taudel]
 gives $$ \phi_n(\alpha)=\prod_{k= 1}^n\frac{(S(a_ke^{-k\tau_n}))^{b_k}}{(S(a_ke^{-k\delta_n}))^{b_k}}
\sim \frac{ \nicef(\tau_n)}{\nicef(\delta_n)}=$$$$ \exp\left(\sum_{l=0}^r h_l\big(\tau_n^{-\rho_l}-\delta_n^{-\rho_l}\big)
-A_0\log\big(\tau_n -\delta_n)+ \Delta(\tau_n)-\Delta(\delta_n)\right),$$\be\ n\to \infty.\la{fkscr}\eu
To bound  $\vert\phi_n(\alpha)\vert,\ \alpha\in [\alpha_0(n),1/2]$ from above, we use the formula
\be \log \vert\phi_n(\alpha)\vert=\Re \big(\log\phi_n(\alpha)\big).\la{relog}\eu
We see that \ $\text{for\ \ all }\ \ \alpha\in[\alpha_0(n),1/2),$
 $$\Re \big(\tau_n^{-\rho_l}\big)=\Re\Big((\delta_n-2\pi  \alpha \ i)^{-\rho_l}\Big)= \Re\Big(\delta_n^{-\rho_l}(1- \frac{2\pi  \alpha \ i}{\delta_n})^{-\rho_l}\Big)=$$
$$
\delta_n^{-\rho_l}\left(\Big(1+\big(\frac{2\pi \alpha}{\delta_n}\big)^2\Big)^{\frac{-\rho_l}{2}}\Re\big(w_n(\alpha)\big)\right)\le  $$\be \delta_n^{-\rho_l}\Big(1+ (2\pi)^2\delta_n^{\rho_r}\log^2 n\Big)^{\frac{-\rho_l}{2}},
l=1,2,\ldots r, \la{sju}\eu
where $w_n(\alpha)$ is a complex variable with
$\vert w_n(\alpha)\vert=1.$

Continuing the last inequality we have, for $n$ sufficiently large and $\alpha\in [\alpha_0(n),1/2]$:
\be \Re\Big((\delta_n-2\pi  \alpha \ i)^{-\rho_l}\Big)\le \delta_n^{-\rho_l}\big(1-
\frac{\rho_l}{2}(2\pi)^2\delta_n^{\rho_r}\log^2 n\big),\ l=1,\ldots,r.\eu
Consequently, for all $\alpha\in[\alpha_0(n),1/2)$ and  $l=1,\ldots,r-1$,
$$ \Re \big(
\tau_n^{-\rho_l}\big)-\delta_n^{-\rho_l}\le -C\delta_n^{\rho_r-\rho_l}\log^2 n\to 0,\ n\to \infty, \ C=\frac{\rho_1}{2}(2\pi)^2>0,$$
while \be     \Re \big(\tau_n^{-\rho_r}\big)-\delta_n^{-\rho_r}\le -C\log^2 n,\ n\to \infty,\ C=\frac{\rho_r}{2}(2\pi)^2>0,\la{vertlog}\eu
with equality for $\alpha=\alpha_0(n)$. Also, we have
\be -A_0\Re\big(\log \frac{\tau_n}{\delta_n}\big)= -\frac{A_0}{2}\log\Big(1+\big(\frac{2\pi\alpha}{\delta_n}\big)^2\Big).\la{gucs}\eu
Finally, it follows from \refm[fkscr], \refm[vertlog] and \refm[gucs] that for all $\alpha\in[\alpha_0(n),1/2),$
\be \Big\vert \log\vert \phi_n(\alpha)\vert \Big\vert\ge C\log^2 n,\ n\to \infty,\ C>0 \la{kbdcuk}\eu
and
  $$\vert \phi_n(\alpha)\vert\le\big(1+(\frac{2\pi \alpha}{\delta_n})^2\big)^{-\frac{A_0}{2}}\exp\big(-C\log^2 n))=o((\delta_n)^{1+\frac{\rho_r}{2}}),$$
 \be n\to \infty,\ C>0.\la{juac}\eu
  \refm[juac] yields \refm[relat] and, consequently  \refm[locl], which proves  that  \refm[nicef1] implies  \refm[ygb]. The sufficiency of \refm[ygb] for LLT is immediate.
\qed

$(iv)$ Proof of the conditions \refm[,nur] and \refm[nscon].

 Our  first goal  is to bound from above the function $\Big\vert\log\vert\phi_n(\alpha)\vert\Big\vert$ for \ rational $\alpha\in [\alpha_0(n),1/2].$
By virtue of \refm[relog] and \refm[ch12] we write \be
\Big\vert\log\vert\phi_n(\alpha)\vert\Big\vert=\frac{1}{2}\sum_{k=1}^n b_k\log \Big(\frac{S^2(a_ke^{-k\delta_n})}{\vert S(a_ke^{-k\delta_n+ 2\pi i\alpha k})\vert^2}\Big),\ \alpha\in {\cal R}.\la{work}\eu
We firstly indicate  the following two essential facts:
\be (a)\ U_n(k;\alpha):=\log\Big(\frac{S^2(a_ke^{-k\delta_n})}{\vert S(a_ke^{-k\delta_n+ 2\pi i\alpha k})\vert^2}\Big)
\to 0,\ \la{ebh}\eu\be \text{for\ all}\ \alpha\in {\cal R}\ \text{and\ for \ all}\ k=k(n): k(n)\delta_n\ge C\log^2 n, \ C>0,\ n\to \infty,\ \text{since}\ S(0)=1 \nonumber\eu
 and
\be (b)\ U_n(k;\alpha)=0,\ \text{if} \ \alpha k \ \text{is \ an \ integer}.\ \ \ \ \ \ \ \ \ \ \ \ \ \ \ \ \ \ \ \ \
\eu
As a particular case of $(a)$,
 \be U_n(k;\alpha)\to 0,\ n\to \infty,\ \alpha\in {\cal R}, \text{\ for\ \ all } k\ge \delta_n^{-1-\epsilon},\ \ \text{with \ any }\ \epsilon>0.\la{zbhj}\eu
Next we show that for $n$ sufficiently large, the main contribution to the sum in \refm[work] comes from the terms with  $k\in\kappa_v,\ v=1,\ldots, q-1,$ where, given an integer $q>1$ and $\epsilon>0$, the set  of integers $\kappa_v=\kappa_{v,q}(\epsilon)$ is defined by
$$\kappa_v
:=\{1\le k\le \delta_n^{-1-\epsilon}, \epsilon>0: \ k\delta_n<\infty,\ \text{as}\ n\to \infty,\  \text{and}\ k\equiv v(mod\  q) \}.$$
Let  $0<\alpha\le 1/2$ be a rational number, i.e. $\alpha=\frac{p}{q}>0, \ gcd(p,q)=1,\ q>1.$
For  $n$ large enough, we can assume that $\frac{p}{q}\in [\alpha_0(n),1/2],$ by the above definition of $\alpha_0(n)$. Then
\be U_n(k;\frac{p}{q})\to \log\Big(\frac{S^2(c_k)}{\vert S(c_ke^{ 2\pi i\frac{v}{q}})\vert^2}\Big),\ k\in\kappa_v,\ n\to \infty,\la{unk}\eu

where   $c_k:=\limsup_{n\to \infty}\big( a_ke^{-\delta_n k}\big),\ \text{for } \  k\in \kappa_v.$

Due to our assumption \refm[aks] on $a_k,\ k\ge 1,$ and the definition of the set $\kappa_v,$ the constants $0 \le c_k\le 1, \ k\in \kappa_v$.
 Also,
\be U_n(k;\frac{p}{q})\ge 0,\ k\in \kappa_v, \la{eur}\eu
with equality if and only if in the expansion \refm[szs],\ $d_j=0,\ \text{for \ all}\ j:q\!\!\not| j,$
for some integer $q\ge 2$.
This says that in the case of equality, the condition \refm[,nur] of Theorem 2 does not hold, which leads to
 $\vert \phi_n(\frac{p}{q}) \vert= 1,$ in contradiction to  the condition \refm[ygb]. Hence, the  condition \refm[,nur] is necessary for $LLT$ to hold, which means that   \refm[,nur] guarantees the strict inequality in \refm[eur].\\
We now prove the necessity of the bounds \refm[nscon].
 Let $k=ql+v\in \kappa_v: \delta_n k= c+ \epsilon_n,$ with $c\ge 0,\ \epsilon_n\to 0, \ as\ n\to\infty,$ and let $z_0= e^{-c}a_k e^{ 2\pi i\frac{v}{q}}$ be a zero of some  order $m\ge 1$ of  the function $S(z),\ z=e^{-\delta_n k}a_k  e^{ 2\pi i\frac{v}{q}}:$
$$S(z)=(z-z_0)^m\tilde{ S}(z),\ \tilde{ S}(z_0)\neq 0. $$

Recalling (see the discussion after \refm[logSexp]) that it should be $\vert z_0\vert=1,$ it follows that in the above representation of $z_0$, $c=0, a_k=1,$
and therefore, in the case discussed,
$$\vert S(e^{-\delta_n k}  e^{ 2\pi i\frac{v}{q}})\vert\sim(\delta_n k)^m \vert\tilde{S}(z_0)\vert= O((\delta_nk)^m)\le C\delta_n,\ C>0,$$\be k\in \kappa_v: \delta_n k=  \epsilon_n\to 0,\ n\to \infty. \la{xupv}\eu

Finally, if $z_0\neq 1$ is a pole of $S(z)$ in the circle $\vert z\vert\le 1$, then again it should be $\vert z_0\vert=1,$ by virtue of the aforementioned remark, and therefore, $c=0,a_k=1.$ Consequently, applying \refm[rgul] with $l=1,$ we get
$\vert S(e^{-\delta_n k}e^{2\pi i\frac{v}{q}})\vert\sim C(\delta_n k)^{-1},\ k\in \kappa_v: \delta_nk\to 0, n\to \infty$, for some $C>0$.
Taking into account that $\vert S(z_0)\vert\le S(1),$ we conclude that in the case considered $z=1$ should be a pole as well, which implies
$$ U_n(k;\frac{p}{q})\sim\log\frac{\vert L^2(1)\vert(\delta_nk)^{-2}}{\vert L^2(e^{2\pi i\frac{v}{q}})\vert(\delta_nk)^{-2}}=\log\frac{\vert L(1)\vert}{\vert L(e^{2\pi i\frac{v}{q}})\vert}=O(1),\ n\to \infty.$$

In  light of the above reasoning  it is clear that the following two cases should be broadly distinguished, as  $n\to \infty$:
  $$ (A): 
 0< U_n(k;\frac{p}{q})\le C>0$$
 \be (B): 
0<U_n(k;\frac{p}{q})
\le -C\log \delta_n,\ C>0,\la{n,fho}\eu
where $\frac{p}{q}\in [\alpha_0(n),1/2]$ and where in both cases
the bounds  hold  for all integers $q\ge 2$ and for all $$k\in \kappa(q;\epsilon):=\{1\le k\le \delta_n^{-1-\epsilon}: q\!\!\not| k, k\delta_n<\infty,\ n\to \infty \}.$$

 Recalling that $S(z)$ may have only a finite number of zeros and poles on the unit disk, it follows that the above bounds on $U_n$ are valid for all poles and zeros with the same constant.

{\bf Remark} The case $(A)$ is in force when $S(z)$ has a complex pole on the boundary of the unit disk, while the case $(B)$ holds when $S(z)$ has a complex zero on the boundary of the unit disk.

Combining \refm[kbdcuk],\refm[zbhj],\refm[n,fho] and \refm[work] we  derive the required bound:
\be \Big\vert\log\vert\phi_n(\frac{p}{q})\vert\Big\vert\le
    \big(\sum_{1\le k\le n,q\!\!\not| k} b_k \big)\left\{
                                                    \begin{array}{ll}
                                                      C, & \hbox{in Case A} \\
                                                      -C\log \delta_n, & \hbox{in\ Case B}
                                                    \end{array}
                                                  \right.
\la{gucsv}
  \eu
  as $n\to \infty$.\ \refm[gucsv] together with the lower bound
\refm[kbdcuk] proves the necessity of  conditions \refm[,nur] and \refm[nscon] for $NLLT$.

Note that the inequality \refm[gucsv] holds for all $\alpha\in [\alpha_0,1/2),$ because  the characteristic function $\phi_n(\alpha)$ is continuous in $\alpha\in {\cal R}$ for any given $n\ge 1.$

The sufficiency of \refm[nscon] and \refm[,nur] for $NLLT$ follows from the fact that the bound \refm[kbdcuk] is necessary and sufficient for \refm[relat].

{\bf Remark}  Condition \refm[nscon] was motivated by Example 2 in \cite{GSE}, which says
that if  the sequence of weights $\{b_k,\ k\ge 1\}$ is supported on the set of $k=lq,\ l=1,2,\ldots,$ with some integer $q>1,$ then $NLLT$ does not hold, though Conditions $I$ and $II$ may hold. In this regard, the condition \refm[nscon] determines the minimal "total mass" of weights $b_k$ that should be concentrated on integers $k\le n$ that are not divisible by a given $q>1,$
 in order that  $NLLT$  be in  force.

\begin{itemize}

\item {\bf Corollary}
\end{itemize}

The proof is simple. Substituting into  \refm[rep] the relation \refm[asnicef], the formula \refm[nicef1] with
 $\delta=\delta_n$ and the asymptotic formula \refm[locl] for $NLLT,$ gives the asymptotic formula \refm[expgr].

\vskip .25cm
{\bf IV. Examples }
\vskip .5cm

{\bf Example 1: Partitions into primes.} This example demonstrates situation when the asymptotic expansion \refm[niceff]   of $\log~\nicef(\delta),\ \delta\to 0$, is not possible.
For partitions into primes,
 $$f(z)=\prod_{k\ge 1}(1-z^k)^{-b_k},\ \vert z\vert<1, $$ where $$ b_k=\left\{
                                                                                                         \begin{array}{ll}
                                                                                                           1, & \hbox{if }\ k \ is\ a\ prime \\
                                                                                                           0, & \hbox{otherwise.}
                                                                                                         \end{array}                                                                                                       \right.
$$
Correspondingly, (see \cite{vaugh},p.117),
\be D_b(s)= \sum_{p} \frac{1}{p^s}:= P(s),   \eu
where the summation is over the set of all primes, is the so- called Prime zeta function that admits the representation
\be P(s)=\sum_{k\ge 1}\mu(k)\frac{\log\zeta(ks)}{k}, \la{prime}\eu
where $\mu$ is the M$\ddot{o}$bius function. It is clear from \refm[prime] that $P(s)$ is analytic in the half-plane $\Re(s)>1$ and that it has a meromorhic continuation to the strip $0<\Re(s)<1$. Also, it follows from \refm[prime] that $P(s)$ has  logarithmic singularities at the following points: $(i)$  $s=1$, where
$$P(s)\sim \log\zeta(s)\sim \log\frac{1}{s-1}, \ s\to 1,$$  $(ii)$  $s=\frac{1/2+i\rho}{k},\ k\ge 1$ which are induced by   non-trivial zeros $\rho$ of $\zeta(s)$ and $(iii)$  $s=\frac{1}{k},$ with square-free integers $k$, which are induced by those  trivial zeros $-k$ of $\zeta(s)$ at which $\mu(k)\neq 0$.

As a result, the line $\Re(s)=0$ is the natural barrier for the Dirichlet series for $P(s)$, which means that the series cannot be continued analytically to the left of the line $\Re(s)=0$. The latter together
with the fact that the singularities of $P(s)$ are not poles prevents the application of Meinardus' approach.
Roth and Szekeres \cite{szek} were able to obtain the principal term of the asymptotics of $\log c_n, \ n\to \infty,$ with the help of a complicated analysis adapted for the model considered
(see also \cite{vaugh}).

{\bf Example 2:  Partitions into powers of primes.}

Though the structure  in the example below is not exponential, because the condition $I$ of Theorem 1 does not hold, it is possible to apply the technique of Meinardus to obtain the main term in the expansion
of $\log\nicef(\delta)$.
The generating function for the model considered is
\be f(z)=\prod_{k\ge 1}(1-z^k)^{-b_k},\ \vert z\vert<1, \la{prim} \eu

\be b_k=\left\{
                                                                                                         \begin{array}{ll}
                                                                                                           log\ p, & \hbox{if }\ k=p^r,\ \text{where}\ p\text\ {is \ a\ prime\ and\ r\ is\ an \ integer }
  \\
                                                                                                           0, & \hbox{otherwise.}
                                                                                                         \end{array}
\la{afr}                                                                                                       \right.
\eu
The weighted  generating function in \refm[afr] was suggested in 1950 by Brigham (for references see \cite{yif}).

It is known that $b_k=\Lambda(k),$ where $\Lambda(k),\ k\ge 1$ is the von Mangoldt function,  the Dirichlet generating
function of which is \be D_b(s)=-\frac{\zeta^{\prime}(s)}{\zeta(s)}\la{mang}.\eu

Thus, in the case considered $D(s)=\zeta(s+1)D_b(s),$ where  $D_b(s)$ is a meromorhic function in ${\cal C}$ having  poles at all trivial  and non-trivial zeros of $\zeta(s).$
Since the non-trivial zeros are known to be complex numbers, which location depends on the solution of the Riemann hypothesis,
the condition $I$ of our Theorem 1 does not hold.
 However, we will show that the technique of the present paper applied to
 the function $Q(s):=\delta^{-s} \Gamma(s)D(s)$ allows to find the main term in the asymptotic expansion for the $\log~\nicef(\delta),$ as $\delta\to 0^+,$ recovering the results of Richmond \cite{richm}, and Yang \cite{yif}.

$(i)$ By the functional equation for $\zeta$- function,
\be \zeta(s)=2^s\pi^{s-1}\sin(\frac{\pi s}{2})\Gamma(1-s)\zeta(1-s),\ s\in {\cal C}.\eu
Expanding $\sin(\frac{\pi s}{2})$ around the  trivial zeros $\{-2k,\ k\ge 1\}$ of $\zeta(s)$
and taking into account that the function
$2^s\pi^{s-1}\Gamma(1-s)\zeta(1-s)$ is analytic at these points and does not equal to $0$, shows that the function $\frac{1}{\zeta(s)}$ has simple poles at $s=-2k,\ k\ge 1$. Consequently, the function $Q(s)$ defined above, has at each of the  above points a pole of the second order, with residue $O(\delta^{2k} \log \delta)\to 0, \delta\to 0^+ $, so that these poles influence the remainder term $\Delta(\delta)$ only, and the same  is also true for the simple poles $\{-2k-1,\ k\ge 0\}$ with the residues $O(\delta^{2k+1})\to 0,\ \delta\to 0^+,\ k=0,1,2,\ldots,$ induced by  $\Gamma(s)$;

$(ii)$ Recalling the Laurent series expansion for the Riemann zeta function:  $$ \zeta(s)=\frac{1}{s-1}+\sum_{l=0}^\infty \frac{(-1)^l}{l!} \gamma_l \; (s-1)^l,$$ where $\gamma_0=\gamma$ is the Euler constant and $\gamma_l,\ l=1,2,\ldots$ are the Stiltjes constants, we have
from \refm[mang]:
\be D_b(s)=\frac{1}{s-1}+ O(s-1),\ s\to 1.\eu
This shows that $D_b(s)$ has also a simple pole at $s=1,$ with the residue $1$. As a result, the residue of $Q(s)$ at $s=1$ is
equal to $\delta^{-1}\zeta(2)$;

$(iii)$ $s=0$ is the second order pole of $Q(s)$ with the residue $$-\frac{\zeta^\prime(0)}{\zeta(0)} \log \delta + const;$$

$(iv)$ Non-trivial zeros of $\zeta(s).$  These zeros are known to belong  to  the critical strip  $0\le \Re(s)< 1$.
 We adopt the argument in \cite{yif}, Section 5, which is  based on the preceding works of Richmond, Brigham, Ingham and other researchers.
The key observation is that denoting by $\rho$ the non-trivial zeros of $\zeta(s)$, the sum of residues of $Q(s)$ at these points is equal to
\be \sum_{\rho}\delta^{-\rho} \Gamma(\rho)\zeta(1+\rho),\la{tpx}\eu
where the sum is taken over all $\rho$ counted with their multiplicities.
Next, using the known bound for the number of the non-trivial zeros $\rho: Im(\rho)\in [T,T+1] $ one gets that the sum
in \refm[tpx] is $\lll\delta ^{-\theta}$ for some $1/2\le \theta<1$.

Summarizing  $(i)$ -$(iv)$ it follows that the rightmost pole in the model is $s=1$, so that $\delta_n\sim (\zeta(2))^{1/2}n^{-\frac{1}{2}},$  by \refm[asdeltan]. Consequently, by the same argument as in \refm[intrep3],\refm[tsdr]  we deduce from \refm[tsdr] the asymptotic formula by Richmond (see Theorem B from \cite{yif}):
\be \log~\nicef(\delta_n)= 2(\zeta(2))^{1/2}n^{\frac{1}{2}}+ O(n^{\theta/2}),\ 1/2\le\theta<1,\ n\to \infty,  \eu
where $\theta=1/2,$ if the Riemann Hypothesis is true.

Finally, applying  the known asymptotic relation for the von Mangoldt function $\sum_{k=1}^x \Lambda(k)\sim x,\ x\to \infty$ shows that the condition \refm[nscon] for the weights $b_k=\Lambda(k)$ holds, so that NLLT is in force.

{\bf Example 3}
In the example below we build a model that satisfies conditions $I$ and $II$, but disobeys   condition \refm[nscon] of Theorem 2.
So, for this model the generating function grows exponentially, while the local limit theorem does not hold. The idea of the construction is  motivated by  Example 3 in \cite{GSE}.
Let $S(z)=(1-z)^{-1},\ \vert z\vert<1,\ \ a_k\equiv 1,\ k\ge 1$ and
\be
b_k=\left\{
      \begin{array}{ll}
        (k\log^{\epsilon}k)^{-1},\ 0<\epsilon <1, & \hbox{if}\ 4\!\!\not| k,\ k\ge 2,\\
        (k\log^{\epsilon}k)^{-1}+1,\ 0<\epsilon < 1, & \hbox{if} \ 4 | k.
      \end{array}
    \right.
\eu
Thus, in the case considered, the function $D_b(s)$ can be written as
$$D_b(s)= D_b(s;\epsilon)=D_b^{(1)}
(s;\epsilon)+ 4^{-s}\zeta(s),$$
where
$$D_b^{(1)}(s;\epsilon):=\sum_{k=2}^{\infty}\frac{1}{k^{s+1}\log^{\epsilon}k}.$$

Denoting $f(x;s,\epsilon)=\big(x^{s+1}\log^\epsilon x\big)^{-1},$ we apply the  Euler-Maclaurin summation formula (see e.g.\cite{AAR}) to get:

$$ D_b^{(1)}(s;\epsilon)= \int_2^\infty f(x;s,\epsilon)dx + \frac{f(2;s,\epsilon) + f(\infty;s,\epsilon)}{2} +$$\be \int_2^\infty f^{\prime}_x(x;s,\epsilon)(x-[x]) dx,\la{eum}\eu
where $[x]$ is the integer part of $x$. In the  representation \refm[eum], $f(\infty;s,\epsilon)=0,\ \text{for}\ \Re(s)>-1,\ 0<\epsilon<1,$
and $f^{\prime}_x(x;s,\epsilon)=-x^{-s-2}\log^{-\epsilon} x \big(s+1 + \epsilon \log^{-1}x\big).$  The latter implies that the second  integral in \refm[eum] converges absolutely for $\Re(s)> -1 ,\ \epsilon>0.$
Since the first integral in \refm[eum], which we denote $Q(s)$ converges for $s>0$, it is left to show that $Q(s)$ admits meromorhic continuation to all ${\cal C}.$
We have for $\Re(s)>0,\ 0<\epsilon<1,$
$$Q(s):=\int_2^\infty f(x;s,\epsilon)dx = -s^{-1}\big(x^{-s}\log^{-\epsilon}x\big)\vert_2^{\infty}-\epsilon s^{-1}\int_2^{\infty} x^{-s-1}\log^{-\epsilon-1}x\ dx=$$$$ s^{-1}2^{-s}\log^{-\epsilon}2- \epsilon s^{-1}\int_2^{\infty} x^{-s-1}\log^{-\epsilon-1}x\ dx,$$ which can be rewritten as a differential equation with respect to $Q(s):$
\be(sQ(s))^{\prime}=-2^{-s}\log^{1-\epsilon}2+ \epsilon Q(s).\la{dif}\eu
Here we made use of the fact that
$\big(\int_2^{\infty} x^{-s-1}\log^{-\epsilon-1}x\ dx\big)^{\prime}_s=-Q(s).$

The solution of \refm[dif] can be found explicitly:
\be Q(s)=s^{\epsilon-1}\Big(C+ (\log^{1-\epsilon}2) \int_0^s2^{-u} u^{-\epsilon}du\Big),\ 0<\epsilon<1,\ \Re(s)\neq 0, \la{tbkh}\eu
where $C=C(\epsilon)$ is a constant.
 We conclude from \refm[tbkh] that $Q(s)$ allows analytic continuation to all $s\in{\cal C}/{0},$
 with residue $0$ at $s=0,$  because
$$\lim_{s\to 0^+} sQ(s)= 0,\quad  0<\epsilon<1.$$
  In view of \refm[eum],  the residue of $D_b^1(s;\epsilon), \ 0<\epsilon<1$ at $s=0$   is  zero, as well. Thus, conditions $I,II$ hold, while the NLLT is not in force, since the condition \refm[nscon] is violated:
$$\sum_{1\le k\le \delta_n^{-1}: \ 4\!\!\not| k} b_k=\sum_{1\le k\le \delta_n^{-1}: \ 4\!\!\not| k}     (k\log^{\epsilon}k)^{-1}=O(-\log^{1-\epsilon}\delta_n)
,\ n\to \infty,\ 0<\epsilon<1.$$

{\bf Remark on models with non exponential rate of growh.} In Applications, mainly in   number theory, one meets  models for which $c_n$ does not grow exponentially with $n$.
Such cases are known long ago, two typical examples are the number of square-free integers and Goldbach partitions of an even integer into sum of two odd primes. For recent developments in the study of such models see \cite{csin} and \cite{mut}.
Note that the second model is an multiplicative one, while the first is not.

\end{document}